\documentclass{article}
\usepackage[margin=1.0in]{geometry}

\author{John K. Sikora \\ john.sikora@xtera.com}
\title{Analysis of the High Water Mark Convergents of Champernowne's Constant in Various Bases}

\usepackage{hyperref}   
\usepackage{xspace}     
\usepackage{amsthm}     
\usepackage{amsmath}    
\usepackage{fixltx2e}   
\usepackage{twoopt}     
\usepackage{verbatim}   
\usepackage{mathtools}  
\usepackage[table]{xcolor} 
\usepackage{longtable}  
\usepackage{multirow}   
\usepackage{hhline}     
\usepackage{cleveref}   

\newcommand{\cc}{Champernowne's Constant\xspace}
\newcommand{\ct}{\ifmmode C_{10} \else $C_{\mathit{10}}$\xspace \fi}      
\newcommand{\cb}{\ifmmode C_b \else $C_b$\xspace \fi}  
\newcommand{\po}{position\xspace} 
 
\newcommand{\cn}[1][]{coefficient number#1\xspace} 
\newcommand{\Cn}[1][]{Coefficient number#1\xspace} 
\newcommand{\mni}[1]{\mathit{#1}} 
\newcommandtwoopt{\tbp}[2][b][X]{\ifmmode \mni{10}_{\mni{#1}}^{\mni{#2}} \else $\mni{10}_{\mni{#1}}^{\mni{#2}}$\xspace \fi} 
\newcommandtwoopt{\ptbx}[2][b][X]{\ifmmode P\mni{10}_{\mni{#1},\mni{#2}} \else $P\mni{10}_{\mni{#1},\mni{#2}}$\xspace \fi}
\newcommandtwoopt{\dbx}[2][b][X]{\ifmmode D_{\mni{#1},\mni{#2}} \else $D_{\mni{#1},\mni{#2}}$\xspace \fi}
\newcommandtwoopt{\nbx}[2][b][X]{\ifmmode N_{\mni{#1},\mni{#2}} \else $N_{\mni{#1},\mni{#2}}$\xspace \fi}
\newcommand{\csb}[1]{\ifmmode C_{\mni{#1}} \else $C_{\mni{#1}}$\xspace \fi}
\newcommand{\snd}{2\textsuperscript{nd}\xspace}
\newcommand{\sghwm}{2\textsuperscript{nd} Generation HWM\xspace}
\newcommand{\sghwms}{2\textsuperscript{nd} Generation HWMs\xspace}

\let\oldFootnote\footnote
\newcommand\nextToken\relax

\renewcommand\footnote[1]{%
    \oldFootnote{#1}\futurelet\nextToken\isFootnote}

\newcommand\isFootnote{%
    \ifx\footnote\nextToken\textsuperscript{,}\fi}

\newtheoremstyle{jsconj}
{18pt}
{8pt}
{\itshape}    
{}
{\bfseries}
{:}
{1em}
{}

\newtheoremstyle{jsdefn}
{10pt}
{10pt}
{}
{}
{\bfseries}
{:}
{1.5em}
{}

\newtheoremstyle{jsexception}
{6pt}
{6pt}
{}
{}
{\itshape}
{}
{1em}
{}

\theoremstyle{jsconj}
\newtheorem{conj}{Conjecture}
\theoremstyle{jsdefn}
\newtheorem{defn}{Definition}   
\theoremstyle{jsexception}
\newtheorem{excep}{Exception}[conj]
\theoremstyle{plain}

\begin{document}
\numberwithin{equation}{section}

\maketitle 

\begin{abstract}
\noindent In this paper, we show that patterns exist in the properties of the High Water Mark (HWM) convergents of \cc in various bases (\cb), specifically in bases 2 through 124. The convergents are formed by truncating the Continued Fraction Expansion (CFE) of \cb immediately before the CFE HWMs. These patterns have been extended from the known patterns in the CFE HWMs of \cc in base ten (\ct). We show that the patterns may be used to efficiently calculate the CFE coefficients of \cb, and to calculate and predict the lengths of the HWM coefficients, the number of correct digits of \cb as calculated by the convergent, and the convergent error. We have discovered that minor corrections to the pattern formulations in base 10 are required for some bases, and these corrections are presented and discussed. The resulting formulations may be used to make the calculations for \cb in any base.    
\end{abstract}

\section{Introduction} \label{sec:intro}
\cc, \ct, was formulated in 1933 by British mathematician and economist D. G. Champernowne as an example of a normal number \cite{bbcn}. It is formed by concatenating the positive integers to the right of a decimal point, i.e., $0.123456789101112 \ldots$, without end. In 1937, Kurt Mahler proved that \ct is transcendental \cite{bbcnt}.

It is well known that the CFE of \ct consists mostly of coefficients with a reasonably small number of digits, interspersed with coefficients with a very large number of digits. A coefficient that has a value that is larger than any previous coefficient is called a HWM. As one progresses along the CFE, the pattern of relatively small coefficients sprinkled with very large coefficients continues between the HWMs, even though the other large coefficients are not themselves HWMs.

There are many patterns in the properties of the convergents formed by truncating the CFE of \ct immediately before the HWMs for HWM numbers greater than 3. These properties and their patterns have been formulated and discussed \cite[Conjectures~1-7]{bbjks}. The purpose of this paper is to extend the formulations and the conjectures from base 10 to bases 2 through 124, with exceptions noted and the pertinent formulae provided. 

In addition, another conjecture which describes the lengths of the largest coefficients between the HWMs \cite[Conjecture~8]{bbjks2} is generalized from \ct to \cb. We show that deviations from the expected behavior occur in some bases and these deviations are noted and briefly discussed.

\section{Definitions} \label{sec:term}
{\setlength{\parindent}{0cm}

\begin{defn} \label{def:dcc}    
\boldmath \ct \unboldmath - \cc in base 10.
\end{defn}

\begin{defn} \label{def:dcb}
\boldmath \cb \unboldmath - \cc in base $b$. For example \csb{3} is $0.12101112202122100101 \ldots_3$
\end{defn}

\begin{defn} \label{def:dpo}    
\boldmath \cb \unboldmath \textbf{\po}, or \textbf{\po} - For an integer, it is the \po of the first digit of the integer in the consecutive sequence used to generate \cb, starting at \po 0 for the 0 to the left of the radix point. The radix point is not counted. Therefore, for \ct, the integer 1 is in \po 1 and the integer 12 is in \po 14. For $C_3$, the integer $11_3$ is in \po 5.
\end{defn}

\begin{defn} \label{def:dcn}    
\textbf{CFE \cn}, or \textbf{\cn} - The number of the CFE coefficient, starting at \cn 0. \Cn 0 is always 0 for \cb.
\end{defn}

\begin{defn} \label{def:dhwmn}    
\textbf{HWM number}, or \textbf{HWM \#} - The number of the HWM. The first HWM is 0 and it is HWM number 1.
\end{defn}

\begin{defn} \label{def:dhwml}    
\textbf{HWM length} - The number of digits of the HWM of \cb \emph{as represented in base} $b$.
\end{defn}

\begin{defn} \label{def:dtbp}    
\boldmath \tbp \unboldmath - This is equal to $b^X$ and it is read one-zero, or ``ten'', to the power of $X$ in base $b$. This notation is used instead of $b^X$ primarily because the computer programs that were used in the analysis calculate the HWM convergents with two primary inputs: base $b$, and the power of ``ten'', $X$. This terminology was chosen so that there would be consistency between this paper and the programs. 
\end{defn}

\begin{defn} \label{def:dpotbp}    
\boldmath \ptbx \unboldmath \hspace{-0.4pc}, or \boldmath \cb \unboldmath \textbf{\po of} \boldmath \tbp \unboldmath - The formula for the \cb \po of the integer $\tbp, \ X \geq 0$, is:
\begin{equation} \label{eq:cbp}
\cb \text{ \po of } \tbp=\ptbx=1+\sum_{x=0}^X (b-1)\cdot x \cdot b^{(x-1)}
\end{equation}
\end{defn}

Thus, the \cb \po of \tbp[b][0] is 1, the \po of \tbp[b][1] is $b$, and the \po of \tbp[3][2] is 15.

\begin{defn} \label{def:dtbpc}    
\boldmath \tbp \textbf{HWM convergent}, \unboldmath or \boldmath{\tbp} \textbf{convergent} \unboldmath - The convergent formed from the numerator, \nbx, and the denominator, \dbx, as calculated using \Cref{conj:numer} and \Cref{conj:denom}, respectively. The convergent is named in this manner because the \cb \po \tbp is used in the calculation of the convergent. The CFE derived from the convergent does not calculate the \tbp HWM coefficient as defined in \Cref{def:dtbpcc}. The HWM coefficient is, however, the next term in the CFE.
\end{defn} 

\begin{defn} \label{def:dtbpcc}    
\boldmath \tbp \textbf{HWM coefficient}, \unboldmath or \boldmath \tbp \textbf{coefficient} \unboldmath - The CFE coefficient immediately following the last coefficient as calculated by the \tbp convergent. With a few exceptions as noted later in this paper, these coefficients constitute the HWMs of the CFE of \cb.
\end{defn} 

\begin{defn} \label{def:dbx}    
\textbf{The denominator of the} \boldmath \tbp \unboldmath \textbf{HWM convergent}, or \boldmath \dbx \unboldmath  - The denominator for the \tbp HWM convergent as constructed using \Cref{conj:denom}.
\end{defn}

\begin{defn} \label{def:nbx}    
\textbf{The numerator of the} \boldmath \tbp \unboldmath \textbf{HWM convergent}, or \boldmath \nbx \unboldmath  - The numerator for the \tbp HWM convergent as calculated using \Cref{conj:numer}.
\end{defn}

\begin{defn} \label{def:dncd}    
\textbf{Number of correct digits}, or \textbf{NCD}, or $\boldsymbol{NCD(b,X)}$ - The number of correct consecutive digits of \cb as calculated by a convergent, starting at, and counting the 0 to the left of the radix point. The radix point itself is not counted. Therefore, if the \po of the last correct digit, $p_c$, is known, the number of correct digits is $p_c+1$ since the 0 is counted as the first correct digit, but it is in \po 0. When given with arguments, base $b$, and power of ``ten'' $X$, $NCD(b,X)$ is the NCD of the \tbp HWM convergent as defined in \Cref{def:dtbpc}.
\end{defn}
}

In general, in this paper we refer to the HWM convergents in a manner different from that in our previous paper which detailed the HWM convergent properties of \ct \cite{bbjks}. In the previous paper, a HWM convergent is referred to as ``the convergent before HWM \#N''. In this paper a HWM convergent is referred to as ``the \tbp HWM convergent'', or as the ``\tbp convergent'' as defined in \Cref{def:dtbpc}. The reason is that in this paper, bases 2 through 124 are discussed, rather than base 10 exclusively, and the HWM numbers are different for the various bases.\footnote{The HWM numbers match for bases 6 through 124, but not for bases 2 through 5. See \Cref{conj:cncalc} and \Cref{conj:hwml} for details.}

\section{Method} \label{sec:meth}
A program was written in Ruby to calculate \cb in the desired base by taking the digits of each of the integers used to form \cb, dividing by the appropriate power of the base for each digit, and summing. The Euclidean Algorithm applied to real numbers was then used to calculate the CFE coefficients of \cb. The use of this method was difficult and care was taken to calculate \cb to sufficient accuracy so that the desired number of coefficients could be accurately calculated. Bases 2, 7, 8, 9, and 10 were used in the initial investigation.

It was noticed that large coefficients (HWMs) in the CFE of \cb appeared in other bases, similar to the behavior of the CFE of \ct. The CFE coefficients of \cb, terminated before the HWMs were then used to calculate the numerator and denominator of each convergent. Patterns were observed in the denominators of the convergents as given in \Cref{conj:denom}. The means of calculating the numerators of the convergents was determined as given in \Cref{conj:numer}. 

Programs were then written in C to calculate the numerators and denominators using these conjectures and the resulting calculations led to the formulation of \Cref{conj:cncalc} through \Cref{conj:hwml}. These calculations were checked for each \tbp convergent and HWM as described in \Cref{ssec:rv}. 

\section{Conjectures} \label{sec:ce}
In this section we extend the \ct conjectures \cite[Conjectures~1-7]{bbjks} to \cb. The order of the conjectures in the present paper matches the order of the conjectures in the previous paper for ease of comparison. However, \Cref{conj:cncalc} through \Cref{conj:hwml} follow from analyzing the convergents formed using \Cref{conj:denom} and \Cref{conj:numer}.

\begin{conj} \label{conj:cncalc}
The \tbp convergent, for $X \geq 0$, calculates the value of \cb correctly until the integer $\tbp[b][X+1]-\mni{2}$ in the integer sequence used to form \cb is reached. Furthermore, this integer is calculated as  $\tbp[b][X+1]-\mni{1}$ instead of $\tbp[b][X+1]-\mni{2}$. 
\end{conj}

\begin{excep} \label{excep:b2cbc}
\narrower
For base 2, the convergents that calculate \csb{2} per \Cref{conj:cncalc} start with the $\tbp[2][X]$ convergent, $X \geq 2$, rather than $X \geq 0$.  
\end{excep}

\begin{excep} \label{excep:b34cbc}
\narrower
For base 3 and base 4, the HWM convergents that calculate \csb{3} and \csb{4} per \Cref{conj:cncalc} start with the $\tbp[3][X]$ and $\tbp[4][X]$ convergent, $X \geq 1$, rather than $X \geq 0$.  
\end{excep}

\noindent The following are examples of the conjecture: 

\noindent the \tbp[10][1] HWM convergent calculates \ct as ``$0.123456\ldots95969799\ldots$'' instead of \\ ``$0.123456\ldots95969798\ldots$''. Similarly, the \tbp[5][2] HWM convergent calculates \csb{5} as \\ ``$0.12341011\ldots440441442444\ldots_5$'' instead of ``$0.12341011\ldots440441442443\ldots_5$''.

\begin{conj} \label{conj:ncd}
The NCD of the \tbp HWM convergent for base $b, \ X \geq 0$ is given by:
\begin{equation} \label{eq:ncd}
NCD(b,X) = \cb \text{ \po of } \tbp[b][X+1] - X - 2
\end{equation}
\end{conj}
\begin{excep} \label{excep:b2ncd}
\narrower
For $b=2$, \Cref{eq:ncd} is valid for $X \geq 2$.
\end{excep}

\begin{excep} \label{excep:b34ncd}
\narrower
For $b=3$ and $b=4$, \Cref{eq:ncd} is valid for $X \geq 1$.
\end{excep}

\noindent The \cb \po of \tbp[b][X] is given in \Cref{eq:cbp}. Thus, the NCD(10,0) is 8, the NCD(5,3) is 2340, the NCD(16,5) is 99544809, and the NCD(10,7) is 788888881.

\begin{conj} \label{conj:ceven}
The \tbp HWM \cn{s} are even. The \cn{s} start at 0 per \Cref{def:dcn}.
\end{conj}

The \tbp HWM convergents calculate \cb slightly larger than \cb. Per \Cref{def:dcn}, this means that the \tbp HWM \cn{s} will be even. 

This conjecture is very useful at the end of the calculation of the CFE coefficients as the algorithm cannot distinguish the degenerate case where the last coefficient is 1 from the case where the last coefficient is not 1. In the degenerate case where the last coefficient equals 1, the algorithm calculates the last term as $Y+1$ instead of $Y$, and the calculation ends on an even numbered coefficient. However, if the calculation ends on an even numbered coefficient, it is taken from \Cref{conj:ceven} that this is the degenerate case, the even numbered coefficient is $Y$, and the last coefficient (the one before the HWM) is 1.

\begin{conj} \label{conj:error}
The error of the \tbp HWM convergent is given by:
\begin{align} \label{eq:err} 
\text{error} =
\begin{cases}
1.0 \times b^{-(exponent-1)} & \text{base } b \geq 5, X=0 \\
(b-1).1_b \times b^{-exponent} & \text{base } b \geq 3, X=1 \\
1.01_2 \times b^{-exponent} & \text{base } b=2, X=2 \\
(b-1).\underbrace{00\ldots00}_{\mathclap{\text{number of zeroes}}}(b-1)(b-1)_b \times b^{-exponent} & b=2, X \geq 3; \ b \geq 3, X \geq 2 \\
\end{cases}
\end{align}
\begin{align}
\text{number of zeroes} &= X \nonumber \\
exponent &= NCD(b,X)+X+2 \nonumber 
\end{align}
\end{conj}
For example, the \tbp[8][0] HWM has an error of $1.0 \times 8^{-7}$, the \tbp[3][1] HWM has an error of $2.1_3 \times 3^{-15_{10}}$, the \tbp[2][2] HWM has an error of $1.01_2 \times 2^{-18_{10}}$, the \tbp[10][2] HWM has an error of $9.0099_{10}\times 10^{-2890_{10}}$, and the \tbp[15][5] HWM has an error of $\text{e.00000ee}_{15}\times 15^{-67530135_{10}}$. Note that the mantissa is given in base $b$ while the exponent is given in base ten, similar to the output from the GNU Multiple Precision Arithmetic (GMP) library \cite{bbgmp} function \texttt{mpf\_out\_str}. This library was utilized by the programs used in the analysis.

The value of $exponent$ is also equal to the \cb \po of \tbp[b][X+1]. Note that for the base $b=2, X=2$ case we could have explicitly provided the value of $exponent$, but we wanted to show that the value of $exponent$ was consistent for all of the given cases as formulated.

\begin{conj} \label{conj:hwml}
Let:
\begin{align}
(b==2) =
\begin{cases} \nonumber
0 & b \neq 2 \\
1 & b = 2
\end{cases}
\end{align}
\begin{align}
(b==4) =
\begin{cases} \nonumber
0 & b \neq 4 \\
1 & b = 4
\end{cases}
\end{align}

The length, \emph{as represented in base} $b$, of the HWM coefficient that occurs after the last coefficient of the CFE calculated by the \tbp HWM convergent for $X \geq 0$ is given by:
\begin{equation} \label{eq:hwml}
\text{length}=NCD(b,X)-2 \cdot NCD(b,X-1)-3 \cdot (X-(b==2))-2+(b==4)
\end{equation}
\end{conj}

\begin{excep} \label{excep:b2hwml}
\narrower
For base 2, \Cref{eq:hwml} is valid for $X \geq 2$, rather than $X \geq 0$. In addition, the \tbp[2][2] coefficient is \cn 6 and it is of length 3, per \Cref{conj:hwml}. It is not a HWM but the convergent does calculate \csb{2} per \Cref{conj:cncalc}. The \cn of the next HWM, which is the \tbp[2][3] HWM coefficient, is 14.\footnote{The \tbp[2][2] coefficient has a value of 5, one less than the value of \cn[s] 2 and 5, which are the largest terms along with \cn 10 until the \tbp[2][3] HWM coefficient at \cn 14.} All subsequent HWMs are \tbp[2][X] HWM coefficients.
\end{excep}

\begin{excep} \label{excep:b3hwml}
\narrower
For base 3 and base 4, \Cref{eq:hwml} is valid for $X \geq 1$, rather than $X \geq 0$. For base 3, the \tbp[3][1] HWM coefficient is \cn 6, and for base 4, the \tbp[4][1] HWM coefficient is \cn 12. Both are HWMs, and subsequent HWMs are \tbp HWM coefficients.
\end{excep}

\begin{excep} \label{excep:b5hwml}
\narrower
For base 5, the \tbp[5][0] coefficient is \cn 4 and it is of length 1, per \Cref{conj:hwml}.\footnote{The coefficient is 1.} It is not a HWM but the convergent does calculate \csb{5} per \Cref{conj:cncalc}. The \cn of the \tbp[5][1] HWM coefficient is 10. The HWMs at \cn 7 and \cn 9 are not \tbp[5][X] HWM coefficients.
\end{excep}

\noindent For example, the \tbp[2][3] HWM as represented in binary has a length of 9, the \tbp[10][2] HWM as represented in decimal has a length of 2504, and the \tbp[16][5] HWM as represented in hexadecimal has a length of 89198852. 

For the \tbp[b][0] coefficient, $b \geq 5$, the coefficient length is $b-4$. For bases greater than 5, the \tbp[b][0] coefficient is a HWM, and it will be HWM \#4 at \cn 4, per \Cref{def:dcn} and \Cref{def:dhwmn}. All subsequent HWMs are those generated by a \tbp convergent, $X \geq 1$.\footnote{There is a close call for base 6 as \cn 8 is not a \tbp[6][X] coefficient and it is only one less than \cn 4, the \tbp[6][0] coefficient.}  

Furthermore, it appears that for $b \geq 5$, the first 4 CFE terms are [0; $b-2$, $b-1$, 1] = [0; $b-2$, $b$]. This gives a \tbp[b][0] convergent of 
\begin{equation} \label{eq:ptz}
\tbp[b][0] \text{ convergent} =\frac{b}{(b-1)^2}, \ b \geq 5
\end{equation}
As stated above, the \tbp[b][0] convergent is a HWM for $b \geq 6$. For the numerators and denominators for other convergents, see \Cref{eq:numer} and \Cref{eq:denom}, respectively.

It should be noted that for the lowest value of $X$ allowed in any base, the $NCD(b,X-1)$ term is not the true NCD as the \tbp[b][X-1] convergent is not a valid HWM convergent. However, the calculated value as used in \Cref{eq:hwml} gives the correct HWM length. It should also be reiterated that the HWM following the \tbp HWM convergent is not calculated by the convergent. It is calculated by the \tbp[b][X+1] HWM convergent, or by other means \cite{bbdew}. 

It is unknown why the HWM lengths in base 4, as formulated, are longer by 1 digit than the HWM lengths in other bases (other than base 2). It is also unknown why the HWM lengths in base 2, as formulated, are longer by 3 digits than the HWM lengths in other bases (other than base 4). However, in the base 2 case, we feel that since there are no \tbp HWM convergents for $X < 2$, this may increase the lengths of the HWM coefficients. 

\begin{conj} \label{conj:denom}
The denominator, \dbx[b][0], used to calculate the \tbp[b][0] HWM convergent for base, $b \geq 5, \ X=0$, is given in \Cref{eq:ptz}. The denominator used to calculate the \tbp HWM convergent for base, $b, \ X \geq 1$, is given by:
\begin{align} \label{eq:denom}
\dbx =
\begin{cases}
\underbrace{(b-1).(b-1)\ldots (b-1)}_{\mathclap{\text{number of }(b-1)\text{ digits}}}(b-2)\overbrace{0\ldots 0}^{\mathclap{\text{number of zeroes}}}1 \times b^{exponent} & b \text{ is odd} \\
\text{number of } (b-1)\text{ digits} &= X \\
\text{number of zeroes } &= X \\ \\
\left(\frac{b}{2}-1\right).\underbrace{(b-1)\ldots (b-1)}_{\mathclap{\text{number of }(b-1)\text{ digits}}}\overbrace{0\ldots 0}^{\mathclap{\text{number of zeroes}}}\left(\frac{b}{2}\right) \times b^{exponent} & b \text{ is even} \\ 
\text{number of }(b-1)\text{ digits} &= X \\
\text{number of zeroes } &= X+1 \\ \\
exponent &= NCD(b, X-1)+2X+1   
\end{cases}
\end{align}
\end{conj} 

\begin{excep} \label{excep:b2denom}
\narrower
For base 2, \Cref{eq:denom} is valid for $X \geq 2$, rather than $X \geq 1$.\footnote{Although the \tbp[2][1] HWM convergent is not valid, $NCD(2,1)=3$ which gives a valid exponent of 8 for the \tbp[2][2] HWM convergent denominator, $0.110001_2 \times 2^{8_{10}}$.}
\end{excep}

\noindent For example, \dbx[2][4], the denominator for the \tbp[2][4] HWM convergent is $0.1111000001_2 \times 2^{54_{10}}$, \dbx[5][3] is $4.4430001_5 \times 5^{348_{10}}$, \dbx[10][8] is $4.999999990000000005_{10} \times 10^{788888898_{10}}$, and \dbx[15][5] is $\text{e.eeeed000001}_{15} \times \break 15^{3742640_{10}}$.

If $X=1$ for an odd base, then the $(b-1)$ term is to the left of the radix point and the $(b-2)$ term is immediately to the right of the radix point. For example, the denominator for the \tbp[7][1] HWM convergent is $6.501_7 \times 7^{8_{10}}$.

\begin{conj} \label{conj:numer}
The numerator, \nbx[b][0], used to calculate the \tbp[b][0] HWM convergent for base, $b \geq 5, \ X=0$, is given in \Cref{eq:ptz}. For other $b$ and $X$, \nbx is calculated as described below.

\vspace{10pt}
\noindent Given:
\begin{enumerate}
\item The denominator, \dbx, of a \tbp convergent as calculated from \Cref{conj:denom}
\item The \cb \po of the integer \tbp, \ptbx, as given by \Cref{eq:cbp}
\item $C_b[0.1 \ldots \ptbx]$ denotes \cb truncated to \po \ptbx\footnote{The last digit is a 1 as it is the first digit of \tbp in the integer sequence used to form \cb.}
\item $\text{Ceil}(f)$ denotes the integer immediately above non-integer $f$.
\end{enumerate}
Then the numerator, $\nbx, \ X \geq 1$, of the convergent is calculated as:
\begin{align} \label{eq:numer} 
\nbx =
\begin{cases}
\text{Ceil}(\dbx \cdot C_b[0.1 \ldots \ptbx]) + 1 & b \text{ is odd} \\
\text{Ceil}(\dbx \cdot C_b[0.1 \ldots \ptbx]) & b \text{ is even}
\end{cases}
\end{align} 
\end{conj} 

\begin{excep} \label{excep:b2numer}
\narrower
\noindent For base 2, \Cref{eq:numer} is valid for $X \geq 2$, rather than $X \geq 1$.
\end{excep}

For example, the calculation for \nbx[2][3], the numerator for the \tbp[2][3] HWM convergent is: \\
Ceil$(0.11100001_2 \times 2^{21_{10}} \cdot 0.110111001011101111_2) = 110000100000000100001_2$.

The calculation for \nbx[3][2] is: \\
Ceil$(2.21001_3 \times 3^{17_{10}} \cdot 0.121011122021221_3) + 1 = 112222220011222111_3$.

The calculation for \nbx[10][1] is: \\
Ceil$(4.9005_{10} \times 10^{11_{10}} \cdot 0.1234567891_{10}) = 60499999499_{10}$.

And the calculation for \nbx[15][1] is: \\
Ceil$(\text{e.d}01_{15} \times 15^{16_{10}} \cdot 0.123456789\text{abcde}1_{15}) + 1 = \text{120eeeeeeeeeedeed}_{15}$.

\section{Results, Verification, and Data Location} \label{sec:rv}  

\subsection{Results and Verification} \label{ssec:rv}
The \cb CFE coefficients may be calculated very efficiently from the \tbp HWM convergents. For $X \geq 1$, the denominator of the convergent is given by \Cref{conj:denom} and the numerator of the convergent may be calculated by multiplying the denominator by a sufficient number of digits of \cb. The required number of digits is $\ptbx+1$ as given in \Cref{conj:numer}, resulting in a very efficient calculation for the \cb CFE coefficients.\footnote{The required number of digits is $\ptbx+1$ because the zero to the left of the radix point has been included in the count.} The Euclidean Algorithm is then used to calculate the CFE coefficients from the numerator and the denominator. \Cref{conj:ceven} is used to determine if the last coefficient is 1. For $b \geq 5, \ X=0$, \Cref{eq:ptz} is used to calculate the numerator and the denominator.

The properties of the convergents as given in \Cref{conj:cncalc} through \Cref{conj:error} may be observed by dividing the numerator of the convergent by the denominator of the convergent and comparing the result with \cb for each base and power. The property of the convergents as given in \Cref{conj:hwml} may be checked by examination of the \cb CFE coefficients. 

In the current study, the calculations and comparisons were performed for bases 2 to 124. \Cref{tbl:basesx} gives the base, $b$, and the powers of ``ten'', $X$ , in base $b$ that were computed. The calculated results for \textit{each base and each power} were checked against the predictions made in \Cref{conj:cncalc} through \Cref{conj:hwml}. 

The maximum power for each base was usually dictated by memory size. The amount of memory required to calculate \cb, the \cb CFE coefficients, and the error for the \tbp[b][X-1] HWM convergent is approximately the same as the amount of memory required to calculate the \cb CFE coefficients for the \tbp HWM convergent, without calculating \cb and the error. There are separate columns for each type of calculation in \Cref{tbl:basesx}. For more details, see \Cref{ssec:dlads} and \Cref{ssec:c}. 

Sequences A244330 \cite{bboeis2hwmcn}, A244331 \cite{bboeis2hwml} for base 2, and A244332 \cite{bboeis3hwmcn}, and A244333 \cite{bboeis3hwml}, for base 3 have been added to the OEIS which correspond to sequences A038705 \cite{bboeis10hwmcn} and A143534 \cite{bboeis10hwml} in base 10, respectively. The first member of each pair of sequences gives the CFE \cn of each HWM coefficient and the second gives the length of each HWM coefficient as represented in the specified base.\footnote{The \cn{s} start at 1 per the precedent set in A038705 \cite{bboeis10hwmcn}.}  

In addition, sequences A244758 \cite{bboeis2allcl} for base 2, and A244759 \cite{bboeis3allcl} for base 3 have been added to the OEIS which correspond to sequence A143532 \cite{bboeis10allcl} in base 10. These sequences give the number of digits as represented in the specified base for all of the calculated CFE coefficients.

We are considering adding similar sequences for base 4 through base 9. For these bases, there are existing OEIS sequences for the digits of \cb, but not for the CFE coefficients of \cb. Therefore, sequences for the CFE coefficients of \cb would be added as well as sequences for the HWM coefficients, the \cn of the HWMs, the length of all of the calculated coefficients, and the length of the HWM coefficients. 

It should be noted that the maximum numbers of CFE coefficients as calculated and presented in \Cref{tbl:basesx} are not necessarily the maximum number of coefficients that have been calculated in a particular base. For example, at present for \ct, 82328 coefficients have been calculated \cite{bbwmwcncfe}. 

\begin{table}[!htbp]
\centering
\caption{Calculated Bases, Powers, and Number of CFE Coefficients}
\label{tbl:basesx}
\begin{tabular}{|c||c|c|c|c|} 
\hline
\multirow{2}{*}{\textbf{base} $\boldsymbol{b}$} & \multirow{2}{*}{\textbf{Min} $\boldsymbol{X}$} & \multirow{2}{*}{\textbf{Max} $\boldsymbol{X}$ \textbf{w/ error}} & \multirow{2}{*}{\textbf{Max} $\boldsymbol{X}$ \textbf{w/o error}} & \textbf{Number of CFE Coefficients} \\ 
 & & & & \textbf{(for max} $\boldsymbol{X}$ \textbf{w/o error)} \\
\hline
2 & 2 & 25 & 24 & 98093504 \\ \hline  
3 & 1 & 15 & 16 & 2982556 \\ \hline
4 & 1 & 12 & 13 & 629420 \\ \hline
5 & 0 & 10 & 11 & 195554 \\ \hline
6 & 0 & 9 & 10 & 105806 \\ \hline
7 & 0 & 8 & 9 & 53596 \\ \hline

\multirow{2}{*}{8--10} & \multirow{2}{*}{0} & \multirow{2}{*}{7} & \multirow{2}{*}{8} & base 8: 26362 \\
 & & & & base 10: 34062 \\ \hline

\multirow{2}{*}{11--14} & \multirow{2}{*}{0} & \multirow{2}{*}{6} & \multirow{2}{*}{7} & base 11: 14424 \\
 & & & & base 14: 16386 \\ \hline

\multirow{2}{*}{15--16} & \multirow{2}{*}{0} & \multirow{2}{*}{5} & \multirow{2}{*}{6} & base 15: 6080 \\
 & & & & base 16: 6258 \\ \hline

\multirow{2}{*}{17--27} & \multirow{2}{*}{0} & \multirow{2}{*}{4} & \multirow{2}{*}{5} & base 17: 2382 \\
 & & & & base 27: 3184 \\ \hline

\multirow{2}{*}{28--62} & \multirow{2}{*}{0} & \multirow{2}{*}{3} & \multirow{2}{*}{4} & base 28: 1104 \\
 & & & & base 62: 1830 \\ \hline

\multirow{2}{*}{63--90} & \multirow{2}{*}{0} & \multirow{2}{*}{2} & \multirow{2}{*}{3} & base 63: 540 \\
 & & & & base 90: 574 \\ \hline

\multirow{2}{*}{91--124} & \multirow{2}{*}{0} & \multirow{2}{*}{1} & \multirow{2}{*}{2} & base 91: 136 \\
 & & & & base 124: 128 \\ \hline
\end{tabular}
\end{table}

A tabular view of the properties for a specific base, $9$, appears in \Cref{tbl:basenine}. To the best of our knowledge, the values in the shaded cells are predictions, and they have not been confirmed. A similar table, for base 10, is given in \cite[Table 1]{bbjks}.

\begin{table}[!htbp]
\begin{small}
\centering 
\caption{$C_{\mni{9}}$ HWM Convergents Summary}
\label{tbl:basenine}
\begin{tabular}{|c|c|c|c|c|c|c|} 
\hline
 & \textbf{HWM} & \textbf{Integer that} & & & & \\ 
 & \textbf{Coefficient} & \textbf{Fails ; Fails} & & & & \textbf{Length of} \\
 & \textbf{Number} & \textbf{As, in} & & & & \textbf{HWM} \\ 
\textbf{Power of} & \textbf{(next CFE} & \textbf{Convergent} & & \textbf{Convergent} & \textbf{Convergent} & \textbf{(next CFE} \\ 
$\boldsymbol{10_9}$, $\boldsymbol{X}$ & \textbf{term)} & \textbf{Calculation} & $\boldsymbol{NCD(\mni{9},X)}$ & \textbf{Error} & \textbf{Denominator} & \textbf{term)} \\ \hline

\rule{0pt}{2.6ex} 0 & 4 & $7_9;8_9$ & 7 & $1.0 \times 9^{-8_{10}}$ & $64_{10}$ & 5 \\ \hline
\rule{0pt}{2.6ex} 1 & 16 & $87_9;88_9$ & 150 & $8.1_9 \times 9^{-153_{10}}$ & $8.701_9 \times 9^{10_{10}}$ & 131 \\ \hline
\multirow{2}{*}{2} & \multirow{2}{*}{52} & \multirow{2}{*}{$887_9;888_9$} & \multirow{2}{*}{2093} & $8.0088_9$ & $8.87001_9$ & \multirow{2}{*}{1785} \\
 & & & & $\times 9^{-2097_{10}}$ & $\times 9^{155_{10}}$ &  \\ \hline

\multirow{2}{*}{3} & \multirow{2}{*}{152} & \multirow{2}{*}{$8887_9;8888_9$} & \multirow{2}{*}{25420} & $8.00088_9$ & $8.8870001_9$ & \multirow{2}{*}{21223} \\ 
 & & & & $\times 9^{-25425_{10}}$ & $\times 9^{2100_{10}}$ & \\ \hline

\multirow{2}{*}{4} & \multirow{2}{*}{492} & \multirow{2}{*}{$88887_9;88888_9$} & \multirow{2}{*}{287859} & $8.000088_9$ & $8.888700001_9$ & \multirow{2}{*}{237005} \\
 & & & & $\times 9^{-287865_{10}}$ & $\times 9^{25429_{10}}$ & \\ \hline 

\multirow{2}{*}{5} & \multirow{2}{*}{1598} & $888887_9;$ & \multirow{2}{*}{3122210} & $8.0000088_9$ & $8.88887000001_9$ & \multirow{2}{*}{2546475} \\ 
 & & $888888_9$ & & $\times 9^{-3122217_{10}}$ & $\times 9^{287870_{10}}$ & \\ \hline

\multirow{3}{*}{6} & \multirow{3}{*}{4512} & \multirow{2}{*}{$8888887_9;$} & \multirow{3}{*}{32882905} & \multirow{2}{*}{$8.00000088_9$} & $8.888887$ & \multirow{3}{*}{26638465} \\
 &  & \multirow{2}{*}{$8888888_9$} &  & \multirow{2}{*}{$\times 9^{-32882913_{10}}$} & $0000001_9$ &  \\ 
 & & & & & $\times 9^{3122223_{10}}$ & \\ \hline

\multirow{3}{*}{7} & \multirow{3}{*}{12164} & \multirow{2}{*}{$88888887_9;$} & \multirow{3}{*}{338992920} & \multirow{2}{*}{$8.000000088_9$} & $8.8888887$ & \multirow{3}{*}{273227087} \\ 
 &  & \multirow{2}{*}{$88888888_9$} &  & \multirow{2}{*}{$\times 9^{-338992929_{10}}$} & $00000001_9$ &  \\ 
 & & & & & $\times 9^{32882920_{10}}$ & \\ \hline

\multirow{3}{*}{8} & \multirow{3}{*}{30410} & \cellcolor[gray]{0.8} & \cellcolor[gray]{0.8} & \cellcolor[gray]{0.8} & $8.88888887$ & \cellcolor[gray]{0.8} \\ 
 & & \cellcolor[gray]{0.8}\multirow{-2}{*}{$888888887_9;$} & \cellcolor[gray]{0.8} & \cellcolor[gray]{0.8}\multirow{-2}{*}{$8.0000000088_9$} & $000000001_9$ & \cellcolor[gray]{0.8} \\
 & & \cellcolor[gray]{0.8}\multirow{-2}{*}{$888888888_9$} & \cellcolor[gray]{0.8}\multirow{-3}{*}{3438356831} & \cellcolor[gray]{0.8}\multirow{-2}{*}{$\times 9^{-3438356841_{10}}$} & $\times 9^{338992937_{10}}$ & \cellcolor[gray]{0.8}\multirow{-3}{*}{2760370965} \\ \hline
 
 & \tiny{\textbf{\Cref{conj:ceven}}} & \multirow{3}{*}{\tiny{\textbf{\Cref{conj:cncalc}}}} & \multirow{2}{*}{\tiny{\textbf{\Cref{conj:cncalc},}}} & \multirow{2}{*}{\tiny{\textbf{\Cref{conj:cncalc},}}} & \multirow{2}{*}{\tiny{\textbf{\Cref{conj:denom},}}} & \multirow{3}{*}{\tiny{\textbf{\Cref{conj:hwml}}}} \\
 & \tiny{\textbf{(even}} & & \multirow{2}{*}{\tiny{\textbf{\Cref{conj:ncd}}}} & \multirow{2}{*}{\tiny{\textbf{\Cref{conj:error}}}} & \multirow{2}{*}{\tiny{\textbf{\Cref{eq:ptz}}}} & \\
 & \tiny{\textbf{numbered)}} & & & & & \\ \hline 
\end{tabular}
\end{small}
\end{table}

\subsection{Data Location and Digit Symbols} \label{ssec:dlads}
The computed data may be downloaded \cite{bbdataprogloc} and found under the ``data'' folder. The data is available under the Open Data Commons Public Domain Dedication and License (ODC PDDL) \cite{bbodcpddl}.

The downloaded data must be unzipped and extracted. The data includes the numerator of the convergent and the CFE coefficients as calculated from the convergent for each base, $b$, and power of ``ten'', $X$, as given in \Cref{tbl:basesx} under the ``Max $X$ w/o error'' column. For the bases and powers under the ``Max $X$ w/ error'' column, the data also includes \cb as calculated by the convergent as well as the convergent error. However, for $b=2$, $X=25$, the CFE coefficients are not given due to the lengthy calculation time.\footnote{The calculation is in progress.}

The base, $b$, and power, $X$, appear in each file name. The following prefixes identify the types of files:

\vspace{10pt}
\begin{tabular}{llll} 
$\bullet$ & \texttt{cn\_calc\_} & -- & \cb as calculated by the convergent \\
$\bullet$ & \texttt{cn\_cfe\_coeffs\_} & -- & the CFE coefficients calculated from the convergent \\
$\bullet$ & \texttt{cn\_error\_} & -- & \cb as calculated by the convergent minus \cb (approximately 30 digits) \\
$\bullet$ & \texttt{cn\_numer\_} & -- & the numerator of the convergent \\
\end{tabular}
\vspace{10pt}

It should be noted that there are various suffixes, after the base and power, that indicate the specific program that was used to generate the data. These suffixes include, for example, whether or not the error was calculated. 

The symbols used for the digits in various bases are summarized in \Cref{tbl:basessym}. Of course, for any base, $b$, the symbol range is valid only for digits 0 through $b-1$. Referring to the table, it may be seen, for example, that when examining a data file for a base 36 calculation, the symbol for the digit $35_{10}$ is ``z''. However, when examining a data file for a base 37 calculation, the symbol for the digit $35_{10}$ is ``Z'' and the symbol for the digit $36_{10}$ is ``a''.  When examining a  data file for a base 64 calculation, the symbol for the digit $63_{10}$ is ``$|$1'', with the symbol representing a single character. 

\begin{table}[!htbp]
\centering
\caption{Digit Symbols for the Various Bases}
\label{tbl:basessym}
\begin{tabular}{|c||c|c|c|c|c|c|} 
\hline
 & \multicolumn{6}{c|}{\textbf{Digit Symbols}} \\ \hline
\textbf{Base} $\boldsymbol{b}$ & \textbf{0--9} & \textbf{10--35} & \textbf{36--61} & \textbf{62--71} & \textbf{72--97} & \textbf{98--123} \\ \hline
2--10 & 0--9 & & & & & \\ \hline
11--36 & 0--9 & a--z & & & & \\ \hline
37--62 & 0--9 & A--Z & a--z & & & \\ \hline
63--72 & 0--9 & A--Z & a--z & $|$0--$|$9 & & \\ \hline
73--98 & 0--9 & A--Z & a--z & $|$0--$|$9 & $|$a--$|$z & \\ \hline
99-124 & 0--9 & A--Z & a--z & $|$0--$|$9 & $|$A--$|$Z & $|$a--$|$z \\ \hline
\end{tabular}
\end{table}

Discussions of the results of the processing of the data and compilations of the data may be found in \Cref{ssec:oprogdata}. 

\section{\texorpdfstring{\snd}{2nd} Generation HWM Length Extension to Bases 2 through 62}
\label{sec:sg}
It is well known that in addition to the \ct HWM CFE coefficients, there are other coefficients that have a large number of digits that are not themselves HWMs. For example, coefficient 101 (starting from 0) is 140 digits long, but it is not a HWM since \cn 40 is 2504 digits long and \cn 18 is 166 digits long. Coefficient 101 is an example of what we have designated a \sghwm \cite{bbjks2}.

It is known that for \ct, the lengths of the \sghwms appear to be related to the lengths of the HWMs \cite[Conjecture 8]{bbjks2}. The differences in the lengths of a HWM and its corresponding \sghwm show a regular pattern, increasing by 10 for each HWM number increase. We have found that, in general, the same holds true for bases 2 through 62: the difference in the lengths of a HWM and its corresponding \sghwm increases by 10 for each HWM number increase.\footnote{The length is the number of digits as represented in the base of \cb.}\footnote{No analysis was performed for bases 63 through 124.} However, there were irregularities in some bases as can be seen by the shaded cells in \Cref{tbl:bdhwml}, which is a compilation of the data. Note, for example, the entries for bases 4, 5, 9, 22, and 36. 

It should also be noted that starting at base 24, the first \sghwm starts after HWM \#5 rather than after HWM \#6. Further analysis may show that this is true for base 23 and possibly base 22, but this would most likely involve looking at the error of the potential \sghwm convergents \cite[Conjecture 9]{bbjks2}. 

\begin{center}
\begin{longtable}{|c||c|c|c|c|c|c|c|c|c|c|c|c|c|c|} 
\caption{Bases, Difference Between HWM Lengths and \sghwm Lengths} \label{tbl:bdhwml} \\ \hline
\textbf{base} & \multicolumn{14}{c|}{\textbf{HWM Length - \sghwm Length}} \\ \hline \endfirsthead
\multicolumn{15}{c}{\tablename\ \thetable{} -- continued from previous page} \\ \hline 
\textbf{base} & \multicolumn{14}{c|}{\textbf{HWM Length - \sghwm Length}} \\ \hline 
\endhead
\multicolumn{15}{c}{continued on next page} \\
\endfoot 
\endlastfoot   
2 & 80 & 90 & 100 & 110 & 120 & 130 & 140 & 150 & 160 & 170 & 180 & 190 & 200 & 210 \\ \hline
3 & 48 & 58 & 68 & 78 & 88 & 98 & 108 & 118 & 128 & 138 & 148 &  &  & \\ \hline
4 & \cellcolor[gray]{0.8}36 & 47 & 57 & 67 & 77 & 87 & 97 & 107 & 117 &  &  &  &  & \\ \hline
5 & \cellcolor[gray]{0.8}25 & 37 & 47 & 57 & 67 & 77 & 87 & 97 &  &  &  &  &  & \\ \hline
6 & 25 & 35 & 45 & 55 & 65 & 75 & 85 &  &  &  &  &  &  & \\ \hline
7 & 27 & 37 & 47 & 57 & 67 & 77 &  &  &  &  &  &  &  & \\ \hline
8 & 25 & 35 & \cellcolor[gray]{0.8}42 & 55 & 65 & 75 &  &  &  &  &  &  &  & \\ \hline
9 & \cellcolor[gray]{0.8}24 & 37 & 47 & 57 & 67 &  &  &  &  &  &  &  &  & \\ \hline
10 & 26 & 36 & 46 & 56 & 66 &  &  &  &  &  &  &  &  & \\ \hline
11 & 27 & 37 & 47 & 57 & 67 &  &  &  &  &  &  &  &  & \\ \hline
12 & 26 & 36 & 46 & 56 &  &  &  &  &  &  &  &  &  & \\ \hline
13 & 27 & 37 & 47 & 57 &  &  &  &  &  &  &  &  &  & \\ \hline
14 & 26 & 36 & 46 & 56 &  &  &  &  &  &  &  &  &  & \\ \hline
15 & 27 & 37 & \cellcolor[gray]{0.8}44 & 57 &  &  &  &  &  &  &  &  &  & \\ \hline
16 & 26 & 36 & 46 & 56 &  &  &  &  &  &  &  &  &  & \\ \hline
17 & 27 & 37 & 47 &  &  &  &  &  &  &  &  &  &  & \\ \hline
18 & \cellcolor[gray]{0.8}24 & 36 & 46 &  &  &  &  &  &  &  &  &  &  & \\ \hline
19 & 27 & 37 & 47 &  &  &  &  &  &  &  &  &  &  & \\ \hline
20 & 26 & 36 & 46 &  &  &  &  &  &  &  &  &  &  & \\ \hline
21 & 27 & 37 & 47 &  &  &  &  &  &  &  &  &  &  & \\ \hline
22 & 26 & 36 & \cellcolor[gray]{0.8}43 &  &  &  &  &  &  &  &  &  &  & \\ \hline
23 & 27 & 37 & 47 &  &  &  &  &  &  &  &  &  &  & \\ \hline
24 & 16 & 26 & 36 & 46 &  &  &  &  &  &  &  &  &  & \\ \hline
25 & 17 & 27 & 37 & 47 &  &  &  &  &  &  &  &  &  & \\ \hline
26 & \cellcolor[gray]{0.8}14 & \cellcolor[gray]{0.8}24 & 36 & 46 &  &  &  &  &  &  &  &  &  & \\ \hline
27 & 17 & 27 & 37 & 47 &  &  &  &  &  &  &  &  &  & \\ \hline
28 & 16 & 26 & 36 &  &  &  &  &  &  &  &  &  &  & \\ \hline
29 & \cellcolor[gray]{0.8}15 & 27 & 37 &  &  &  &  &  &  &  &  &  &  & \\ \hline
30 & 16 & 26 & 36 &  &  &  &  &  &  &  &  &  &  & \\ \hline
31 & 17 & 27 & 37 &  &  &  &  &  &  &  &  &  &  & \\ \hline
32 & 16 & 26 & \cellcolor[gray]{0.8}34 &  &  &  &  &  &  &  &  &  &  & \\ \hline
33 & 17 & 27 & 37 &  &  &  &  &  &  &  &  &  &  & \\ \hline
34 & 16 & 26 & 36 &  &  &  &  &  &  &  &  &  &  & \\ \hline
35 & 17 & \cellcolor[gray]{0.8}25 & 37 &  &  &  &  &  &  &  &  &  &  & \\ \hline
36 & \cellcolor[gray]{0.8}15 & 26 & 36 &  &  &  &  &  &  &  &  &  &  & \\ \hline
37 & 17 & 27 & \cellcolor[gray]{0.8}34 &  &  &  &  &  &  &  &  &  &  & \\ \hline
38 & 16 & 26 & 36 &  &  &  &  &  &  &  &  &  &  & \\ \hline
39 & 17 & 27 & 37 &  &  &  &  &  &  &  &  &  &  & \\ \hline
40 & 16 & 26 & 36 &  &  &  &  &  &  &  &  &  &  & \\ \hline
41 & 17 & 27 & 37 &  &  &  &  &  &  &  &  &  &  & \\ \hline
42 & 16 & 26 & 36 &  &  &  &  &  &  &  &  &  &  & \\ \hline
43 & \cellcolor[gray]{0.8}14 & \cellcolor[gray]{0.8}25 & 37 &  &  &  &  &  &  &  &  &  &  & \\ \hline
44 & 16 & 26 & 36 &  &  &  &  &  &  &  &  &  &  & \\ \hline
45 & 17 & 27 & 37 &  &  &  &  &  &  &  &  &  &  & \\ \hline
46 & 16 & 26 & 36 &  &  &  &  &  &  &  &  &  &  & \\ \hline
47 & 17 & 27 & 37 &  &  &  &  &  &  &  &  &  &  & \\ \hline
48 & 16 & 26 & 36 &  &  &  &  &  &  &  &  &  &  & \\ \hline
49 & 17 & 27 & 37 &  &  &  &  &  &  &  &  &  &  & \\ \hline
50 & \cellcolor[gray]{0.8}15 & 26 & 36 &  &  &  &  &  &  &  &  &  &  & \\ \hline
51 & 17 & 27 & 37 &  &  &  &  &  &  &  &  &  &  & \\ \hline
52 & 16 & \cellcolor[gray]{0.8}24 & 36 &  &  &  &  &  &  &  &  &  &  & \\ \hline
53 & 17 & 27 & 37 &  &  &  &  &  &  &  &  &  &  & \\ \hline
54 & 16 & 26 & 36 &  &  &  &  &  &  &  &  &  &  & \\ \hline
55 & 17 & 27 & 37 &  &  &  &  &  &  &  &  &  &  & \\ \hline
56 & 16 & 26 & 36 &  &  &  &  &  &  &  &  &  &  & \\ \hline
57 & \cellcolor[gray]{0.8}16 & 27 & 37 &  &  &  &  &  &  &  &  &  &  & \\ \hline
58 & 16 & 26 & 36 &  &  &  &  &  &  &  &  &  &  & \\ \hline
59 & 17 & 27 & 37 &  &  &  &  &  &  &  &  &  &  & \\ \hline
60 & \cellcolor[gray]{0.8}15 & \cellcolor[gray]{0.8}24 & 36 &  &  &  &  &  &  &  &  &  &  & \\ \hline
61 & 17 & 27 & 37 &  &  &  &  &  &  &  &  &  &  & \\ \hline
62 & 16 & 26 & 36 &  &  &  &  &  &  &  &  &  &  & \\ \hline
\end{longtable}
\end{center}

\section{Notes on Computing} \label{sec:nc}
\subsection{Ruby Programs} \label{ssec:r}
As mentioned in \Cref{sec:meth}, a program was written in Ruby to look for the patterns in the CFE coefficients of \cb in various bases. In particular, the HWM coefficients were noted and another Ruby program was written to calculate the numerators and the denominators, as represented in the particular base, of the convergents terminated immediately before the HWM coefficients. The patterns in the denominators and the methods of calculating the numerators as given in \Cref{conj:denom} and \Cref{conj:numer}, respectively, were determined. The main CFE calculating programs were then written in C due to speed considerations. Other programs were written in Ruby and C to compile the statistics from the C program outputs and to formulate and check the conjectures.

\subsection{C Programs} \label{ssec:c}
The main \cb CFE coefficient calculation programs were written in C. Both the MPIR 32 bit \cite{bbmpir} and GMP 64 bit multiple precision libraries were used. In general, the GMP library was used for the highest power of ``ten'', $X$, for each base due to the speed advantage of GMP on Linux versus MPIR on Windows and also due to the larger data handling capability of the 64 bit GMP library versus the 32 bit MPIR library.

There are three main types of programs. The first two main types of programs each have two subtypes since the limit on the input and output base range of the MPIR and GMP libraries is 2 to 62, while the study was conducted to base 124. The last main type finds use only for low bases and large values of $X$, so the extension to higher bases is unnecessary. The following list shows the main types and subtypes of programs:

\begin{itemize}
\item Programs that calculate the convergent, the CFE coefficients of the convergent, \cb as calculated by the convergent, and the convergent error
  \begin{itemize}
	\item Programs for bases 2 through 62
	\item Programs for bases 63 through 124 
	\end{itemize}
\item Programs that calculate the convergent and the CFE coefficients of the convergent
  \begin{itemize}
	\item Programs for bases 2 through 62
	\item Programs for bases 63 through 124
	\end{itemize}
\item A program that calculates the convergent, \cb as calculated by the convergent, and the convergent error
\end{itemize}

\Cref{tbl:basest} gives the calculation times for the various bases and powers of ``ten'', $X$. See \Cref{ssec:hw} for more details on the calculations. It may be observed that the base 63 through base 124 calculation times increase rapidly with $X$ as the base conversion processing is not as efficient as the MPIR and GMP input and output functions. 

Currently, for bases from approximately 4 to 62, the limitation in these calculations for the maximum power, $X$, for a given base, $b$, is usually memory. However, for the lower bases, particularly base 2 and base 3, and higher powers, the number of coefficients that are calculated becomes very high. In these cases, the calculation time of the CFE coefficients becomes a major factor. As more memory becomes available, the CFE coefficient calculation time may become more of a factor for bases higher than base 3 at the higher powers.

The program that calculates the convergent, \cb as calculated by the convergent, and the convergent error, but not the CFE coefficients was written to further check the validity of the applicable conjectures for the highest power of `` ten'' in the lower bases without the very time consuming calculation of the coefficients. This program may also be used for higher bases but the calculation time savings is not nearly as beneficial.

For bases 63 to 124, the limitation is calculation time rather than memory due to the inefficient base conversion. However, as the base increases, the difference in the number of digits between powers of ``ten'' increases, and this factors into the increase in calculation time as well.

\begin{table}[!htbp]
\centering
\caption{Calculation Times}
\label{tbl:basest}
\begin{tabular}{|c||c|c|c|c|} 
\hline
 & & \multicolumn{3}{c|}{\textbf{calculation time in seconds}} \\ \hline
\multirow{3}{*}{\textbf{base} $\boldsymbol{b}$} & \multirow{3}{*}{$\boldsymbol{X}$} & \textbf{CFE} & \textbf{CFE} & \textbf{With} \\ 
 & & \textbf{Calculation} & \textbf{Calculation} & \textbf{Error, No} \\
 & & \textbf{with Error} & \textbf{without Error} & \textbf{Coefficients} \\
\hline
\multirow{2}{*}{2} & 25 &  &  & 219 \\  \cline{2-5}  
 & 24 &  & 2041980 & 105 \\ \hline
\multirow{2}{*}{3} & 15 & 21551 &  & \\ \cline{2-4}
 & 16 &  & 176523 & \\ \hline
\multirow{2}{*}{4} & 12 & 6048 &  & \\ \cline{2-4}
 & 13 &  & 59031 & \\ \hline 
\multirow{2}{*}{5} & 10 & 1206 &  & \\ \cline{2-4}
 & 11 &  & 13169 & \\ \hline 
\multirow{2}{*}{6} & 9 & 869 &  & \\ \cline{2-4}
 & 10 &  & 8908 & \\ \hline 
\multirow{2}{*}{7} & 8 & 397 &  & \\ \cline{2-4}
 & 9 &  & 3161 & \\ \hline 

\multirow{4}{*}{8--10} & \multirow{2}{*}{7} & base 8: 56 &  & \\
 & & base 10: 928 &  & \\ \cline{2-4}
 & \multirow{2}{*}{8} & & base 8: 569 & \\
 & & & base 10: 5524 & \\ \hline

\multirow{4}{*}{11--14} & \multirow{2}{*}{6} & base 11: 130 &  & \\
 & & base 14: 903 &  & \\ \cline{2-4}
 & \multirow{2}{*}{7} & & base 11: 476 & \\
 & & & base 14: 3330 & \\ \hline

\multirow{4}{*}{15--16} & \multirow{2}{*}{5} & base 15: 60 &  & \\
 & & base 16: 36 &  & \\ 
\hhline{~---~}
 & \multirow{2}{*}{6} & & base 15: 139 & \\
 & & & base 16: 152 & \\ \hline

\multirow{4}{*}{17--27} & \multirow{2}{*}{4} & base 17: 5 &  & \\
 & & base 27: 88 &  & \\ 
\hhline{~---~}
 & \multirow{2}{*}{5} & & base 17: 9 & \\
 & & & base 27: 172 & \\ \hline

\multirow{4}{*}{28--62} & \multirow{2}{*}{3} & base 28: 1 &  & \\
 & & base 62: 82 &  & \\ 
\hhline{~---~}
 & \multirow{2}{*}{4} & & base 28: 3 & \\
 & & & base 62: 144 & \\ \hline

\multirow{4}{*}{63--90} & \multirow{2}{*}{2} & base 63: 5271 &  & \\
 & & base 90: 51821 &  & \\ 
\hhline{~---~}
 & \multirow{2}{*}{3} & & base 63: 14275 & \\
 & & & base 90: 53217 & \\ \hline

\multirow{4}{*}{91--124} & \multirow{2}{*}{1} & base 91: 1 &  & \\
 & & base 124: 3 &  & \\ 
\hhline{~---~}
 & \multirow{2}{*}{2} & & base 91: 1 & \\
 & & & base 124: 3 & \\ \hline

\end{tabular}
\end{table}

\subsection{Program Location and Preferred Programs} \label{ssec:plapp}
All of the programs that were used in the study, including ones that are not mentioned, are available under the GNU LGPL \cite{bbgnulgpl}. To download the programs, see reference \cite{bbdataprogloc}, under the ``programs'' folder. As there are many programs, some of which were used to check results, the preferred program for each main type used in the HWM convergent calculations is indicated in \Cref{tbl:pp}. The preferred programs for MPIR are Microsoft Visual Studio C++ projects, while the GMP programs consist of a C program and a C header file. For instructions on how to run the programs, see the \texttt{readme\_mpir.txt} and \texttt{readme\_gmp.txt} files that may be found in the appropriate folders. 

\subsection{Other Programs and Data of Interest} \label{ssec:oprogdata}
As already mentioned, several programs were written to process the data and compare the output of the main programs with the conjectures. One processed data set that we feel may be of general interest is a list of the CFE coefficient lengths for all \cn[s], not just the HWMs. A separate file has been created for each base, $b$, and the maximum power of ``ten'', $X$,  in which the coefficients were calculated. In some of the lower base cases, where the files may be very large (over 1 GB), compilations for lower values of $X$ are included.

Once the \texttt{vcpp.tar.gz} file has been unzipped and extracted, these files may be found in the following folders:

\begin{itemize}
\item \texttt{get\_hwm\_lengths\_generic\_base}
\item \texttt{get\_hwm\_lengths\_generic\_base\_auto}
\item \texttt{get\_hwm\_lengths\_generic\_base\_v\_hi\_base\_auto}
\end{itemize} 

\noindent The name of the file includes $b$ and $X$. An example name is \\ \texttt{cn\_cfe\_coeffs\_base\_2\_pow\_10\_24\_no\_err\_all\_coeffs\_lengths.txt}. 

There are other files in each of these directories for each $b$ and $X$. These files have been processed from the files containing all of the coefficient lengths. The files with \texttt{\_powers} at the end of the name give the HWM lengths for each $b$ and $X$.

\begin{table}[!htbp]
\centering
\caption{Preferred Programs}
\label{tbl:pp}
\begin{tabular}{|c|c|c|} 
\hline
\textbf{Type of Program} & \textbf{Applicable Bases} & \textbf{Preferred Program} \\ \hline
\multirow{4}{*}{w/ error} & \multirow{2}{*}{2--62} & MPIR: \texttt{cn\_generic\_base\_w\_err\_no\_itoa\_auto} \\
 & & GMP: \texttt{cn\_generic\_base\_gmp\_w\_err\_auto} \\ \cline{2-3}
 & \multirow{2}{*}{63--124} & MPIR: \texttt{cn\_generic\_base\_w\_err\_very\_hi\_base\_auto} \\
 & & GMP: \texttt{cn\_generic\_base\_gmp\_w\_err\_very\_hi\_base\_auto} \\ \hline
\multirow{4}{*}{w/o error} & \multirow{2}{*}{2--62} & MPIR: \texttt{cn\_generic\_base\_no\_err\_no\_itoa\_auto} \\
 & & GMP: \texttt{cn\_generic\_base\_gmp\_no\_err} \\ \cline{2-3}
 & \multirow{2}{*}{63--124} & MPIR: \texttt{cn\_generic\_base\_no\_err\_very\_hi\_base\_auto} \\
 & & GMP: \texttt{cn\_generic\_base\_gmp\_no\_err\_very\_hi\_base\_auto} \\ \hline
w/ error, no coeffs & 2--62 & GMP: \texttt{cn\_generic\_base\_gmp\_w\_err\_auto\_no\_coeffs\_calc} \\ \hline
\end{tabular}
\end{table}

\subsection{Hardware} \label{ssec:hw}
The results given in \Cref{tbl:basest} were obtained on a desktop computer with a 3.40 GHz i7$-$2600 CPU, running a single thread, with 8 GB of 1333 MHz DDR3 SDRAM. The programs were run with GMP on Linux, and GMP was built via the command line using \texttt{$-$$-$build=x86\_64$-$intel$-$linux$-$gnu} to optimize GMP to 64 bit operation on an Intel x86 processor. Additional tuning per GMP instructions was attempted, but no performance increase was observed. Building GMP in this manner resulted in code that ran approximately twice as fast as the generic C code.

The MPIR programs were run in 32 bit mode on Windows and the build mode was for generic C code. The GMP programs ran about 8 times faster than the MPIR programs.

\section{Suggestions for Future Analysis} \label{sec:fa}
\subsection{Extending the powers of ``ten''} \label{ssec:ept}
Naturally, extending the maximum power of ``ten'', $X$, for each base and checking the calculations against the conjectures is desirable. Generally, the amount of memory required increases by a factor of approximately $1.25 \cdot base$ for each increment in power. When extending the calculations to a higher power of $X$, we feel that it is very important to make the calculations with the program that calculates \cb as calculated by the convergent as well as with the program that calculates the CFE coefficients without calculating \cb, with no more than one power difference between the two. We believe that ensuring that the calculation of \cb conforms to \Cref{conj:cncalc} for power $X-1$ increases the credibility of the calculations of the CFE coefficients for power $X$.

As explained in \Cref{ssec:c}, for the lower bases, computing time becomes problematic. In general, for bases 62 and lower, each increment in power increases the time by a factor of approximately $2.5 \cdot base$, with a minimum of approximately a factor of 8 for base 2 and high powers of $X$.

\subsection{Analysis of the ratio of adjacent HWM coefficient numbers} \label{ssec:hwmsp}
We feel that it would be interesting to compare the ratio of the \cn{s} of successive HWMs in the various bases. A base dependent formula may exist for the limit of this ratio as the power of ``ten'', $X$ increases without bound.  

\subsection{More efficient algorithm for calculating the CFE coefficients} \label{ssec:cfee}
Implementation of a faster, more efficient algorithm to calculate the CFE coefficients would be of great benefit for the lower bases and higher powers of ``ten''.

\subsection{Further analysis of the \texorpdfstring{\snd}{2nd} Generation HWM coefficients} \label{ssec:sghwma}
Further analysis of the \sghwms is suggested to see if the difference in HWM and \sghwm lengths stabilizes at 10 as the powers of ``ten'' go higher. It would also be of interest to see if the \sghwm convergent error, and denominators in other bases have similar patterns to those in base 10 \cite{bbjks2}.

\subsection{More efficient base conversion for bases greater than base 62} \label{ssec:mebc}
The current programs for bases greater than base 62 have inefficient methods for calculating \cb and for converting between bases. This results in greatly increased calculation times for higher powers of ``ten''. Therefore, it is desired that more efficient algorithms be employed.  

\section{Conclusions} \label{sec:cc}
We have determined that patterns exist in the CFE HWM convergents for \cc in bases 2 to 124, similar to the patterns that exist in the CFE HWM convergents for \cc in base 10. These patterns have been presented in the form of formulae in conjectures. In some bases, slight changes in the formulae are required. There are also a few exceptions to the conjectures, depending upon the base and the power of ``ten''. These differences and exceptions have been presented and discussed. 

A cursory analysis of the \sghwm lengths in base 2 to base 62 has been presented and the results and exceptions to the expected results have been discussed. It was suggested that further study of the \sghwm convergents would be worthwhile.

The programs used in the analysis, and the data calculated by the programs, were discussed. The location of the programs and data was given. The programs are available under the GNU LGPL. The data is available under the ODC PDDL.

\clearpage 

\def\bibindent{1em}

\end{document}